\documentclass[11pt]{amsart}
   \usepackage{amsmath,amssymb}
   
   \renewcommand{\bf}{\bfseries}

   \newcommand{\al}{\alpha}
   
   \newcommand{\T}{\mathbb{T}}
   \newcommand{\Z}{\mathbb{Z}}
   \newcommand{\N}{\mathbb{N}}
   \newcommand{\E}{\mathbb{E}}
   \newcommand{\R}{\mathbb{R}}

   \newtheorem{theorem}{Theorem}
   \newtheorem{lemma}{Lemma}
   \newtheorem{definition}{Definition}
   \newtheorem{proposition}[theorem]{Proposition}
   
   \renewcommand{\epsilon}{\varepsilon}

   \newcommand{\dis}{\displaystyle}
  
\begin{document}
  \title{Averages along cubes for not necessarily commuting measure preserving transformations}
   \author{I. Assani}
\thanks{Department of Mathematics,
UNC Chapel Hill, NC 27599, assani@math.unc.edu}
\thanks{AMS subject classification 37A05, 37A30}

   \begin{abstract}
   We study the pointwise convergence of some weighted averages linked to averages along cubes.
   We show that if $(X,\mathcal{B},\mu, T_i)$ are not necessarily commuting measure
   preserving systems on the same finite measure space and if $f_i,$ $1\leq i\leq 6$ are
    bounded functions then the averages
   $$\frac{1}{N^3}\sum_{n, m, p=1}^N
   f_1(T_1^nx)f_2(T_2^mx)f_3(T_3^px)f_4(T_4^{n+m}x)f_5(T_5^{n+p}x)f_6(T_6^{m+p}x)$$
   converge almost everywhere.
  \end{abstract}
  \maketitle
  \section{Introduction}
   Let $(X,\mathcal{B},\mu, T_i),$ $1\leq i\leq 3,$ be three measure preserving systems on
   the same finite measure space. In \cite{Assani} we proved that if $f_i,$ $1\leq i \leq 3$ are
   three bounded functions then the
   averages $$\frac{1}{N^2}\sum_{n=1}^N
   f_1(T_1^nx)f_2(T_2^mx)f_3(T_3^{n+m}x)$$ converge almost everywhere.  This is
   a bit surprising  as it is known \cite{Berend} that the averages along
   diagonal terms such as $\dis \frac{1}{N}\sum_{n=1}^N
   f_1(T_1^nx)f_2(T_2^nx)$ do not converge even in norm when
   the transformations $T_1$ and $T_2$ do not necessarily commute.
   In the first section of this paper we will extend this result by proving the following theorem.
   \begin{theorem}
   Let $(X,\mathcal{B},\mu, T_i),$ $1\leq i\leq 7,$ be six measure
   preserving systems on the same finite measure space and consider $f_i,$ $1\leq i\leq 6$
    bounded functions. Then the averages
   $$\frac{1}{N^3}\sum_{n, m, p=1}^N
   f_1(T_1^nx)f_2(T_2^mx)f_3(T_3^px)f_4(T_4^{n+m}x)f_5(T_5^{n+p}x)f_6(T_6^{m+p}x)$$
   converge almost everywhere and in norm.
   \end{theorem}

   The method used to prove this theorem is a combination of the following key
   estimates obtained in \cite{Assani} and the ergodic decomposition.
   \begin{lemma}
 Let $a_n$, $b_n$ and $c_n$, $n\in \N$ be three sequences of
  scalars that we assume for simplicity bounded by one. Then for each $N$ positive integer we have
\[
\begin{aligned}
&\big|\frac{1}{N^2}\sum_{m,n=0}^{N-1} a_n. b_{m}. c_{n+m}\big|^2
\\
&\leq
\min\bigg[\sup_t\big|\frac{1}{N}\sum_{m'=1}^{2(N-1)}c_{m'}e^{2\pi
 im't}\big|^2, \sup_t\big|\frac{1}{N}\sum_{n'=1}^Na_{n'}e^{2\pi
in't}\big|^2, \sup_t\big|\frac{1}{N}\sum_{n'=1}^Nb_{n''}e^{2\pi
in''t}\big|^2\bigg]
\end{aligned}
\]
\end{lemma}

\begin{lemma}
Let $$M_N(A_1, A_2,...,A_7)=
\frac{1}{N^3}\sum_{p,n,m=0}^{N-1}a_{1,p}a_{2,n}a_{3,p+n}a_{4,m}a_{5,n+m}a_{6,p+m}a_{7,n+m+p}$$
 the averages of seven bounded (by one) sequences $A_i=(a_{i,n})$, $1\leq i\leq
 7.$ Let us denote by $\mathcal{G}$ the set of couples of integers
 between 1 and 7, $(i, j),$ which are connected by one of the indices $n, m$ or
 $p.$ Then for each $N$ positive integer we have
\[
\begin{aligned}
&\big|M_N(A_1, A_2, ..., A_7)\big|^2 \\
&\leq C \min_{(i,j)\in\mathcal{G}}\bigg[ \max\big[
\frac{1}{N}\sum_{n=0}^{N-1}\sup_t\bigg|\frac{1}{N}\sum_{m=0}^{N-1}a_{i,m}a_{j,n+m}e^{2\pi
imt}\bigg|^2,
\frac{1}{N}\sum_{n=0}^{N-1}\sup_t\bigg|\frac{1}{N}\sum_{m=0}^{2(N-1)}a_{i,m}a_{j,n+m}e^{2\pi
imt}\bigg|^2\big]\bigg].
\end{aligned}
\]

\end{lemma}

 With these lemmas we will derive the pointwise convergence of Wiener-Wintner types
 of averages that will lead to the conclusion stated in Theorem 1. These
 pointwise results extend Wiener-Wintner classical ergodic
 theorem. (see \cite{Assani2}, for instance for several proofs
 of this Wiener Wintner result).
\vskip1ex
 This is done in a first subsection. In a second subsection we will study
 the problem of recurrence to a single set in the case of three
 transformations. We will be able to extend Khintchine's
 recurrence result by studying for any measurable set $A$ with
 positive measure the positivity of the limit
 $$\lim_N \frac{1}{N^2}\sum_{n,m=0}^{N-1}
   \mu\{A\cap T_1^nA\cap T_2^mA\cap T_3^{n+m}A\}
  >0 $$
  when the transformations are not necessarily commuting.
\vskip1ex
  In the second section of the paper we will look at the
  convergence of weighted averages. For a measure preserving transformation $T$ we
  denote by $\mathcal{K}$ the $\sigma$-algebra spanned by the
  eigenfunctions of $T.$
  The method used in \cite{Assani} to prove the pointwise convegence of averages along the cubes
  for the powers of the same
  measure preserving transformation led to the
  following results.
  \begin{lemma}
Let $(X, \mathcal{B}, \mu, T)$ be an ergodic dynamical system and
let $f\in \mathcal{K}^{\perp}.$ Then for $\mu$ a.e. $x$ for all
bounded sequences $a_n$ , $b_n$ , $c_n$,
\begin{enumerate}
\item $\lim_N\dis\frac{1}{N^2}\sum_{n,m=0}^{N-1}a_nb_mf(T^{n+m}x)=
0$, \item $\dis \lim_N\frac{1}{N^2}\sum_{n,m=0}^{N-1}f(T^nx)b_m
c_{n+m}= 0$ and \item $\lim_N\dis
\frac{1}{N^2}\sum_{n,m=0}^{N-1}a_n f(T^mx)c_{n+m}= 0.$
\end{enumerate}
\end{lemma}
and
\begin{proposition}
Let $(X, \mathcal B, \mu, T)$ be an ergodic dynamical system and
let $f\in L^2(\mu).$  Then for $\mu$ a.e. $x$ for all bounded
sequences $a_n$, $b_n$ such that $\dis
\frac{1}{N}\sum_{n=0}^{N-1}a_ne^{2\pi int}$ and $\dis
\frac{1}{N}\sum_{n=0}^{N-1}b_ne^{2\pi int}$ converge for each $t$,
the sequence
$$\frac{1}{N^2}\sum_{n=0}^{N-1}a_n b_m f(T^{n+m}x)$$ converges.
A similar statement holds if one replaces $a_n$ with $f(T^nx)$ and
uses instead $b_m$ and $c_{n+m}$ or if one chooses $b_m= f(T^mx)$
and uses $a_n$ and $c_{n+m}.$
\end{proposition}
 The intriguing aspect of these results is the fact that the set of convergence for $x$
 is independent of the bounded sequences $a_n$, $b_n$ and $c_n$. An illustration of such property can
 be given by taking
 $a_n= (f_1(T_1^nx)),$ $b_m= f_2(T_2^mx)$ with $f_1, f_2\in
 L^{\infty}.$ One obtains immediately the almost everywhere
 convergence of the averages $$\frac{1}{N^2}\sum_{n,m =1}^{N}f_1(T_1^nx)f_2(T_2^mx)f(T^{n+m}x)$$
 if the transformation $T$ is ergodic.
 Other choices for the sequences $a_n$ and $b_n$ are also possible. For instance one could easily take
 $a_n = f(T^{p(n)})(x)$ where $p(x)$ is a real polynomial with
 positive integer coefficients.
 . Such observation seemed to indicate that the almost everywhere convergence of the averages
  along the cubes of a single transformation, namely
 $\frac{1}{N^2}\sum_{n, m=1}^N f(T^nx)g(T^mx) h(T^{n+m}x),$ relies more on the underlying arithmetic structure than
 on its dynamical structure. This is one of the reasons why we
 asked in \cite{Assani} if the assumption of ergodicity made in
 Lemma 3 and the Proposition 2 above was necessary. In this paper
 we will answer in part this question by showing that the ergodicity
 assumption is indeed necessary in Lemma 3. At the present time we
 do not know if Proposition 2 is true without ergodicity
 assumption. With the method used in \cite{Assani} we have the
 following
 \begin{proposition}
Let $(X, \mathcal B, \mu, T)$ be a measure preserving system and
let $f\in L^2(\mu).$ Define the set $\mathcal{WW}_1$ as
$$\mathcal{WW}_1= \big\{ a\in l^{\infty}; \lim_N
\frac{1}{N}\sum_{n=0}^{N-1} a_n e^{2\pi i nt} \text{exists for all
t}\big\}$$. If the set
$$D= \big\{x: \frac{1}{N^2}\sum_{n,m=0}^{N-1} a_nb_mf(T^{n+m}x),
\text{converge for all bounded sequences $(a_n)\in
\mathcal{WW}_1$, $(b_m)\in \mathcal{WW}_1$}\big\}
$$
is measurable then for $\mu$ a.e. $x$ for all bounded sequences
$(a_n)\in \mathcal{WW}_1$, $(b_m)\in \mathcal{WW}_1$ the averages
$$\frac{1}{N^2}\sum_{n,m=0}^{N-1} a_nb_mf(T^{n+m}x)$$
converge.
\end{proposition}
\vskip1ex

The currently open question is the measurability of $D$ that we
will not address in this paper.

\vskip1ex

  It is worth pointing out that if one looks only at the
 norm convergence Proposition 2 is true without ergodicity.

 \vskip1ex

 We will also look at the higher order averages. We denote by $A_i = (a_{n,i})$
  $1\leq i\leq 6,$ six bounded sequences of scalars. We consider the averages
  $$M_N(A_1, A_2,...A_6,f)(x)= \frac{1}{N^3}\sum_{n,
m,p=0}^{N-1}a_{1,p}a_{2,n}a_{3,p+n}a_{4,m}a_{5,n+m}a_{6,
p+m}f(T^{n+m+p}x).$$  In \cite{Assani} we proved that if $f\in
\mathcal{CL}^{\perp}$ then we have a similar result to Lemma 1.
More precisely we have;
 \begin{proposition}
Let $(X, \mathcal{B}, \mu, T)$ be an ergodic dynamical system and
let $f\in CL^{\perp}$. Then for $\mu$ a.e. $x$ for all bounded
sequences $A_i= (a_{i,n})$, $1\leq i\leq 6$ the sequence
$$M_N(A_1, A_2,...A_6,f)(x)= \frac{1}{N^3}\sum_{n,
m,p=0}^{N-1}a_{1,p}a_{2,n}a_{3,p+n}a_{4,m}a_{5,n+m}a_{6,
p+m}f(T^{n+m+p}x)$$ converge to zero.
\end{proposition}

A natural question is to find the precise condition on the
sequences $A_i$ that will give the almost everywhere convergence
of the averages $M_N(A_1, A_2,...A_6,f)(x)$ when $f\in
\mathcal{CL}$. We will show that a condition such as
$\lim_N\frac{1}{N}\sum_{n=1}^N a_{n,i}e^{2\pi int}$ exists for
each $t\in\R$ which is actually necessary and sufficient for the
convergence of the weighted averages
$$\frac{1}{N^2}\sum_{n=0}^{N-1}a_n b_m f(T^{n+m}x)$$ in the universal sense described by Proposition 2 is no longer
sufficient for the convergence of the averages $M_N(A_1,
A_2,...A_6,f)(x).$ At the present time sufficient conditions on
the sequences that would guarantee the almost convergence are not
yet clear to us.

\section{Almost everywhere convergence and recurrence for not necessarily commuting
measure preserving transformations}
\subsection{Proof of Theorem 1}
We recall that if $(X, \mathcal{B}, \mu, T)$ is a measure
preserving dynamical system then the measure $\mu$ can be
disintegrated in a product so that $d\mu = d\mu_c dc$ and $(X,
\mathcal{B}, \mu_c)$ becomes an ergodic dynamical system. This
disintegration allows to lift several results from the ergodic
case to the not necessarily ergodic one. \vskip1ex The proof of
Theorem 1 will be completed after several steps. First we will
need a Wiener Wintner strengthening of Theorem 10 in
\cite{Assani}.
\begin{lemma}
Let $(X, \mathcal{B}, \mu, T_i)$ be three measure preserving
transformations on the same finite measure space. Consider three
bounded functions $f_i$, $1\leq i\leq 3$ then for $\mu$ a.e. $x$
for all $\epsilon_1, \epsilon_2 \in \R$ the averages
$$\frac{1}{N^2}\sum_{n,m =0}^{N-1}
f_1(T_1^nx)f_2(T_2^mx)f_3(T_3^{n+m}x) e^{2\pi in\epsilon_1}e^{2\pi
im\epsilon_2}$$ converge.
\end{lemma}
\begin{proof}
Without loss of generality we can assume that the functions $f_i$
are bounded by one. We use the ergodic decomposition with respect
to $T_3$ to obtain a disintegration of $\mu$, $d\mu= d\mu_{c,
3}dc$ into ergodic components. By the same disintegration and
because of the Wiener Wintner theorem for measure preserving
transformations for $c$ a.e., for $\mu_{c, 3}$ a.e. $y,$ the
averages $$\frac{1}{N}\sum_{n=0}^{N-1} f_1(T_1^ny)e^{2\pi
in\epsilon_1}$$ and  $$\frac{1}{N}\sum_{m=0}^{N-1}
f_2(T_2^my)e^{2\pi im\epsilon_2}$$ converge for all $\epsilon_1,
\epsilon_2 \in \R$ . It is clear that the transformations $T_1$
and $T_2$ may no longer be measure preserving with respect to
$\mu_{c, 3}$ but we are only using here the disintegration of
measurable sets of full measure given by Wiener Wintner ergodic
theorem. Let us consider the Kronecker factor $\mathcal{K}_{3,c}$
of $T_3$ with respect to $(X, \mathcal{B}, \mu_{c, 3})$ and let us
decompose the function $f_3$ into the sum $f_{3,K_c}+ f_{3,
K_{c}^{\perp}}$ where $f_{3, K_c}$ is its projection onto
$\mathcal{K}_{3,c}.$ By Bourgain's uniform Wiener Wintner ergodic
theorem (see \cite{Assani2} for instance for a proof) we have for
$\mu_{3, c}$ a.e. $y$
$$\lim_N\sup_t\big|\frac{1}{N}\sum_{n=0}^{N-1} f_{3, K_c}(T_3^ny)e^{2\pi
int}\big|= 0 .$$ Applying Lemma 1 with $a_n= f_1(T_1^ny)e^{2\pi
in\epsilon_1},$ $b_m= f_2(T_2^my)e^{2\pi im\epsilon_2},$  and $c_k
= f_3(T_3^ky)$ we obtain the estimate
$$\sup_{\epsilon_1, \epsilon_2}\bigg|\frac{1}{N^2}\sum_{n, m= 0}^{N-1}f_1(T_1^ny)f_2(T_2^my)f_{3, K_{c}^{\perp}}(T_{3}^{n+m}y) e^{2\pi in\epsilon_1}e^{2\pi
im\epsilon_2}\bigg|\leq
\sup_t\big|\frac{1}{N}\sum_{k=0}^{N-1}f_{3,
K_{c}^{\perp}}(T_3^{k}y)e^{2\pi ikt}\big|.$$ As a consequence of
the uniform Wiener Wintner theorem we have
$$\lim_N\sup_{\epsilon_1, \epsilon_2}\bigg|\frac{1}{N^2}\sum_{n, m= 0}^{N-1}f_1(T_1^ny)f_2(T_2^my)f_{3, K_{c}^{\perp}}(T_{3}^{n+m}y) e^{2\pi in\epsilon_1}e^{2\pi
im\epsilon_2}\bigg|= 0.$$ The function $f_{3, K_c}$ projects onto
the eigenfunctions of $T_3$ with respect to $\mu_{c, 3}$. If
$e_{j, 3}$ is one of these eigenfunctions with corresponding
eigenvalue $e^{2\pi i\theta_j}$ then we have
$$f_{3, K_c} = \big[\sum_{j=0}^{\infty} \int f_{3, K_c}(y)
\overline{e_{j,3}}(y)d\mu_{c, 3}(y)e_{j, 3}\big].$$ Hence by
linearity and approximation it is enough to consider the case
where $f_{3, K_c}$ is one of the eigenfunctions $e_{j, 3}.$ In
this case $f_{3, K_c}(T_3^{n+m}y) = e^{2\pi i(n+m)\theta_j}e_{j,
3}$ and the averages become
\[
\begin{aligned}
&e_{j, 3}\frac{1}{N^2}\sum_{n,m =0}^{N-1} f_1(T_1^ny)f_2(T_2^my)
e^{2\pi i(n+m)\theta_j}e^{2\pi in\epsilon_1}e^{2\pi
im\epsilon_2}\\
&= e_{j, 3}\frac{1}{N^2}\sum_{n,m =0}^{N-1} f_1(T_1^ny)e^{2\pi
in(\theta_j+ \epsilon_1)}f_2(T_2^my)e^{2\pi im(\theta_j+
\epsilon_2)}
\end{aligned}
\]
The convergence can be derived now by the disintegration, done at
the beginning of the proof, of the sets where the Wiener Wintner
ergodic theorem applied to the functions $f_1$ and $f_2.$

As the set of $x$ for which  for all $\epsilon_1, \epsilon_2 \in
\R$ the averages
$$\frac{1}{N^2}\sum_{n,m =0}^{N-1}
f_1(T_1^nx)f_2(T_2^mx)f_3(T_3^{n+m}x) e^{2\pi in\epsilon_1}e^{2\pi
im\epsilon_2}$$ converge is $\mathcal{B}$ measurable we can
integrate with respect to $\mu_{c,3}$ and $dc$ to show that this
set has full measure.

\end{proof}

\begin{lemma}
Let $(X, \mathcal{B}, \mu, T_i)$ be four measure preserving
transformations on the same finite measure space. Then for all
bounded functions $f_i$, $1\leq i\leq 4,$ the averages
$$\frac{1}{N^3}\sum_{n, m, p =0}^{N-1}
f_1(T_1^nx)f_2(T_2^mx)f_3(T_3^{n+m}x)f_4(T_4^{p}x)$$ converge
$\mu$ a.e. and in norm.
\end{lemma}
\begin{proof}
We can write these averages as
$$\big[\frac{1}{N} \sum_{p=0}^{N-1} f_4(T_4^px)\big]\big[\frac{1}{N^2}\sum_{n,
m=0}^{N-1}f_1(T_1^nx)f_2(T_2^mx)f_3(T_3^{n+m}x)\big].$$ The
conclusion now follows from Birkhoff's pointwise ergodic theorem
and Theorem 10 in \cite{Assani}. The convergence in norm is an
easy consequence of Lebesgue dominated convergence theorem.
\end{proof}
We need now a Wiener Wintner version of Lemma 5.
\begin{lemma}
Let $(X, \mathcal{B}, \mu, T_i)$ be four measure preserving
transformations on the same finite measure space. Then for all
bounded functions $f_i$, $1\leq i\leq 4,$ for $\mu$ a.e. $x$, for
all $\epsilon_1, \epsilon_2\in \R$ the averages
$$\frac{1}{N^3}\sum_{n, m, p=0}^{N-1}f_1(T_1^nx)f_2(T_2^mx)f_3(T_3^{n+m}x)f_4(T_4^{p}x)
e^{2\pi i(n+p)\epsilon_1}e^{2\pi i(m+p)\epsilon_2}$$ converge.
\end{lemma}

\begin{proof}
 We can rewrite the averages as
$$\big[\frac{1}{N} \sum_{p=0}^{N-1} f_4(T_4^px)e^{2\pi ip(\epsilon_1 +\epsilon_2)}\big]\big[\frac{1}{N^2}\sum_{n,
m=0}^{N-1}f_1(T_1^nx)f_2(T_2^mx)f_3(T_3^{n+m}x)e^{2\pi
in\epsilon_1}e^{2\pi im\epsilon_2}\big].$$ The a.e. convergence is
a consequence of the Wiener Wintner ergodic theorem for measure
preserving transformations and Lemma 4.

\end{proof}

\begin{lemma}
Let $(X, \mathcal{B}, \mu, T_i)$ be five measure preserving
transformations on the same finite measure space. Then for all
bounded functions $f_i$, $1\leq i\leq 5,$ for $\mu$ a.e. $x$ the
averages
$$\frac{1}{N^3}\sum_{n, m,
p=0}^{N-1}f_1(T_1^nx)f_2(T_2^mx)f_3(T_3^{n+m}x)f_4(T_4^{p}x)f_5(T_5^{n+p}x)$$
converge.
\end{lemma}
\begin{proof}
We follow the path of the proof of Lemma 4. The set where the
averages converge is $\mathcal{B}$ measurable. We use the ergodic
decomposition of $(X, \mathcal{B}, \mu, T_5)$ into ergodic
components on $(X, \mathcal{B}, \mu_{c,5})$. We disintegrate the
set where the averages
$$\frac{1}{N^3}\sum_{n, m, p=0}^{N-1}f_1(T_1^nx)f_2(T_2^mx)f_3(T_3^{n+m}x)f_4(T_4^{p}x)
e^{2\pi i(n+p)\epsilon}$$ converge for each $\epsilon\in \R.$ We
decompose the function $f_5$ into its projection onto the
corresponding Kronecker factor $f_{5,K_c}$ and
$f_{5,K_{c}^{\perp}}.$ By considering first the case of one
eigenfunction then by approximation and linearity we obtain for
$\mu_{c, 5}$ a.e. y the convergence of the averages
$$\frac{1}{N^3}\sum_{n, m,
p=0}^{N-1}f_1(T_1^ny)f_2(T_2^my)f_3(T_3^{n+m}y)f_4(T_4^{p}y)f_{5,K_c}(T_5^{n+p}y).$$
We can dominate the averages with the function $f_{5,K_c^{\perp}}$
by their absolute value
$$\bigg|\frac{1}{N^3}\sum_{n, m,
p=0}^{N-1}f_1(T_1^ny)f_2(T_2^my)f_3(T_3^{n+m}y)f_4(T_4^{p}y)f_{5,K_c^{\perp}}(T_5^{n+p}y)\bigg|$$
which in turn are bounded by

$$\bigg|\frac{1}{N^3}\sum_{m=0}^{N-1}|f_2(T_2^my)|\sum_{n=0}^{N-1}|f_1(T_1^ny)||f_3(T_3^{n+m}y)|
|\sum_{p=0}^{N-1}f_4(T_4^{p}y)f_{5,K_c^{\perp}}(T_5^{n+p}y)|.$$
Using the fact that the functions are uniformly bounded (by one
without loss of generality) we get the upper bound
$$\frac{1}{N}\sum_{n=0}^{N-1}\big|\frac{1}{N}\sum_{p=0}^{N-1}f_4(T_4^{p}y)f_{5,K_c^{\perp}}(T_5^{n+p}y)\big|.$$
We can apply the remark made after the proof of Lemma 5 in
\cite{Assani} to obtain the bound
$$\sup_t\big|\frac{1}{N}\sum_{k=0}^{N-1}f_{5,K_c^{\perp}}(T_5^{k}y)e^{2\pi
ikt}\big|$$ which converges to zero by the uniform Wiener-Wintner
ergodic theorem. By combining the convergence obtained for
functions $f_{5,K_c}$ and $f_{5,K_{c}^{\perp}}$ we can reach the
$\mu_{5, c}$ a.e. $y$ convergence of the averages
$$\frac{1}{N^3}\sum_{n, m,
p=0}^{N-1}f_1(T_1^ny)f_2(T_2^my)f_3(T_3^{n+m}y)f_4(T_4^{p}y)f_{5,K_c}(T_5^{n+p}y).$$
The convergence $\mu$ a.e. x can be obtained by integration with
respect to $d\mu_{5,c}dc$ .

\end{proof}
It remains to add one more transformation and function, namely
$T_6$ and $f_6.$ The path is quite clear . We start with a Wiener
Wintner version of the Lemma 7.
\begin{lemma}
Let $(X, \mathcal{B}, \mu, T_i)$ be five measure preserving
transformations on the same finite measure space. Then for all
bounded functions $f_i$, $1\leq i\leq 5,$ for $\mu$ a.e. $x$ for
all $t\in \R$ the averages
$$\frac{1}{N^3}\sum_{n, m,
p=0}^{N-1}f_1(T_1^nx)f_2(T_2^mx)f_3(T_3^{n+m}x)f_4(T_4^{p}x)f_5(T_5^{n+p}x)e^{2\pi
i(m+p)t}$$ converge.
\end{lemma}
\begin{proof}
We will use several tools in the proof of the previous Lemmas 7
and 8.  We reconsider the disintegration of the measure $\mu$ into
ergodic components with respect to $(X, \mathcal{B}, \mu_{5, c}).$
We disintegrate the measurable set where the averages
$$\frac{1}{N^3}\sum_{n, m, p=0}^{N-1}f_1(T_1^nx)f_2(T_2^mx)f_3(T_3^{n+m}x)f_4(T_4^{p}x)
e^{2\pi i(n+p)\epsilon_1}e^{2\pi i(m+p)\epsilon_2}$$ converge for
all $\epsilon_1, \epsilon_2\in \R.$ We disintegrate also the
measurable set where the averages
$$\frac{1}{N^3}\sum_{n, m,
p=0}^{N-1}f_1(T_1^nx)f_2(T_2^mx)f_3(T_3^{n+m}x)f_4(T_4^{p}x)f_5(T_5^{n+p}x)e^{2\pi
i(m+p)t}$$ converge for all $t\in \R.$ We decompose the function
$f_5$ into its projection onto the corresponding Kronecker factor
$f_{5,K_c}$ and $f_{5,K_{c}^{\perp}}.$ Again by approximation and
linearity it is enough to look at the case of an eigenfunction
$e_{j,5}$ with eigenvalue $e^{2\pi i\theta_j}.$ The averages in
this case are equal to
$$\frac{1}{N^3}\sum_{n,m,p=0}^{N-1}f_1(T_1^ny)f_2(T_2^my)f_3(T_3^{n+m}y)f_4(T_4^{p}y)e^{2\pi
i(n+p)\theta_j}e^{2\pi i(m+p)t}$$ and converge $\mu_{c,5}$ a.e.
$y$ for all $t$. We are left with the averages related to the
function $f_{5,K_{c}^{\perp}}.$ By observations similar to those
made in Lemma 7 we obtain for each t the upper bound
$$\frac{1}{N}\sum_{n=0}^{N-1}\big|\frac{1}{N}\sum_{p=0}^{N-1}f_4(T_4^{p}y)e^{2\pi ipt}f_{5,K_c^{\perp}}(T_5^{n+p}y)\big|.$$
This last term is dominated by
$$\sup_s\big|\frac{1}{N}\sum_{k=0}^{N-1}f_{5,K_c^{\perp}}(T_5^{k}y)e^{2\pi
iks}\big|$$ and the convergence follows by the uniform Wiener
Wintner ergodic theorem.

\end{proof}

\noindent{\bf End of the proof of Theorem 1}

\vskip1ex

We consider $(X, \mathcal{B}, \mu, T_i)$ six measure preserving
transformations on the same finite measure space and six bounded
functions $f_i$, $1\leq i\leq 6.$ We want to prove that for $\mu$
a.e. x the averages
$$\frac{1}{N^3}\sum_{n, m,
p=0}^{N-1}f_1(T_1^nx)f_2(T_2^mx)f_3(T_3^{n+m}x)f_4(T_4^{p}x)f_5(T_5^{n+p}x)f_6(T_6^{m+p}x)$$
converge.

\vskip1ex We use an ergodic decomposition of the measure $\mu$
into ergodic components for $T_6$. As in the previous lemmas this
reduces the study of the convergence on these components. The
function $f_6$ is decomposed into the sum $f_{6, K_c}$ and $f_{6,
K_c^{\perp}}.$ The convergence $\mu_{6,c}$ a.e. $y$ is obtained by
linearity, approximation and the use of Lemma 8. It remains to
prove the convergence for the averages related to $f_{6,
K_c^{\perp}}.$ We can use Lemma 2 with the sequence $a_{7,k}= 1$
(see also the proof of Lemma 6 in \cite{Assani} to obtain the
following inequalities:
\[
\begin{aligned}
&\bigg|\frac{1}{N^3}\sum_{n, m,
p=0}^{N-1}f_1(T_1^nx)f_2(T_2^mx)f_3(T_3^{n+m}x)f_4(T_4^{p}x)f_5(T_5^{n+p}x)f_{6,
K_c^{\perp}}(T_6^{m+p}x)\bigg| \\
& \leq \frac{1}{N^3}\sum_{n,m=0}^{N-1}\big|\sum_{p=0}^{N-1}
f_4(T_4^py)f_5(T_5^{n+p}y)f_{6,K_c^{\perp}}(T_6^{m+p}y)\big| \\
&\leq C
\sup_s\big|\frac{1}{N}\sum_{k=0}^{N-1}f_{6,K_c^{\perp}}(T_6^{k}y)e^{2\pi
iks}\big|
\end{aligned}
\]
The conclusion of the theorem follows after using the uniform
Wiener Wintner ergodic theorem and integration.

\subsection{An extension of Khintchine recurrence theorem}
\vskip1ex

 Khintchine classical recurrence theorem says that if $A$ is a
set of positive measure, $T$ an invertible measure preserving
system and $\epsilon>0$ the set
$$\{n\in {\Z}: \int {\mathbf{1}}_A.{\mathbf{1}}_A
   \circ T^n d\mu\geq [\int {\mathbf 1}_A d\mu]^2 - \epsilon\}$$
   has bounded gaps.
This recurrence result states that for any measurable set $A$ with
positive measure its images under the iterates of
   $T$ come back and overlap the set with bounded gaps. This is a
   consequence of von Neumann mean ergodic theorem as $$\lim_{N\rightarrow{\infty}}
   \int { \frac{1}{N}}\sum_{n=1}^N{\mathbf{1}}_A.{\mathbf{1}}_A\circ T^n d\mu\geq \mu(A)^2.$$
In this section we study similar recurrence properties with two
and three measure preserving transformations that do not
necessarily commute. We first give a two dimensional extension of
Khintchine's theorem. We can remark that an example given in
\cite{leibman} shows that the averages
$$\frac{1}{N^2}\sum_{n,m=1}^N\mu(A\cap T_1^{-n}A\cap T_2^{-m}A\cap
T_1^{-n}T_2^{-m}A)$$ may diverge if $T_1$ and $T_2$ do not
necessarily commute.

\begin{proposition}
Let $(X, \mathcal{B}, \mu)$ be a probability measure space and
$T_1$, $T_2$ two measure preserving transformations on this
measure space. We denote by $\mathcal{I}_1$ and $\mathcal{I}_2$
the $\sigma$ algebras of the invariant sets for $T_1$ and $T_2$.
Consider $A$ a set of positive measure. Then
$$\lim_N \frac{1}{N^2}\sum_{n, m=1}^N \mu(A\cap T_1^{-n}A\cap
T_2^{-n-m}A)= \int_A\E(\mathbf{1}_A,
\mathcal{I}_1)(x).\E(\mathbf{1}_A, \mathcal{I}_2)(x)d\mu.$$ In
particular
$$\lim_N \frac{1}{N^2}\sum_{n, m=1}^N \mu(A\cap T_1^{-n}A\cap
T_2^{-n-m}A)\geq \mu(A)^4.$$
\end{proposition}
\begin{proof}
The averages
$$\frac{1}{N^2}\sum_{n, m=1}^N \mu(A\cap T_1^{-n}A\cap
T_2^{-n-m}A)$$ are the integrals of the functions
$$\frac{1}{N^2}\sum_{n, m = 1}^N \mathbf{1}_A(x)\mathbf{1}_A
(T_1^n x)\mathbf{1}_A(T_2^{n+m}x)$$ with respect to the measure
$\mu$. As a particular case of Theorem 1 we have the pointwise
convergence of these averages. Thus $$\lim_N\frac{1}{N^2}\sum_{n,
m=1}^N \mu(A\cap T_1^{-n}A\cap T_2^{-n-m}A)$$ exists after
integration. So we just have to prove that
$$\lim_N\frac{1}{N^2}\sum_{n, m=1}^N
\mathbf{1}_A(T_1^nx)\mathbf{1}_A(T_2^{n+m}x)= \E(\mathbf{1}_A,
\mathcal{I}_1)(x).\E(\mathbf{1}_A, \mathcal{I}_2)(x)$$ in $L^2$
norm to conclude. For each $N$ we have
\[\begin{aligned}
&\frac{1}{N^2}\sum_{n, m=1}^N
\mathbf{1}_A(T_1^nx)\mathbf{1}_A(T_2^{n+m}x)\\
&= \frac{1}{N^2}\sum_{n,m=1}^N
\mathbf{1}_A(T_1^nx)\E(\mathbf{1}_A, \mathcal{I}_2)(x) +
\frac{1}{N^2}\sum_{n, m=1}^N
\mathbf{1}_A(T_1^nx)[\mathbf{1}_A(T_2^{n+m}x)- \E(\mathbf{1}_A,
\mathcal{I}_2)(x)]
\end{aligned}
\]
The first term of the last equation converges by Birkhoff's
pointwise ergodic theorem to $\E(\mathbf{1}_A,
\mathcal{I}_1)(x).\E(\mathbf{1}_A, \mathcal{I}_2)(x)$. Noticing
that the function $\E(\mathbf{1}_A, \mathcal{I}_2)(x)$ is $T_2$
invariant we can bound the $L^2$ norm of the second term by
 $$\|\frac{1}{N}\sum_{n=1}^N\big|\frac{1}{N}\sum_{m=1}^N[\mathbf{1}_A\circ
 T_2^m - \E(\mathbf{1}_A, \mathcal{I}_2)]\circ T_2^n\|_2.$$
 This term is less than
 $$\frac{1}{N}\sum_{n=1}^N \|\sum_{n=1}^N\big|\frac{1}{N}\sum_{m=1}^N[\mathbf{1}_A\circ
 T_2^m - \E(\mathbf{1}_A, \mathcal{I}_2)]\|_2$$
 which is equal to
 $$\|\frac{1}{N}\sum_{m=1}^N[\mathbf{1}_A\circ
 T_2^m - \E(\mathbf{1}_A, \mathcal{I}_2)]\|_2$$
 This last term tends to zero by the mean ergodic theorem applied
 to $T_2$. This proves that
 $\dis \lim_N\|\frac{1}{N^2}\sum_{n, m=1}^N
\mathbf{1}_A(T_1^nx)\mathbf{1}_A(T_2^{n+m}x)- \E(\mathbf{1}_A,
\mathcal{I}_1)(x).\E(\mathbf{1}_A, \mathcal{I}_2)(x)\|_2 = 0.$ It
remains to show that $$\int_A \E(\mathbf{1}_A,
\mathcal{I}_1)(x).\E(\mathbf{1}_A, \mathcal{I}_2)(x)d\mu\geq
\mu(A)^4.$$ We have
\[\begin{aligned}
&\int_A \E(\mathbf{1}_A, \mathcal{I}_1)(x).\E(\mathbf{1}_A,
\mathcal{I}_2)(x)d\mu= \int \E\big(\mathbf{1}_A\E(\mathbf{1}_A,
\mathcal{I}_1),\mathcal{I}_2)\big)\mathbf{1}_Ad\mu\\
&= \int \E\big(\mathbf{1}_A\E(\mathbf{1}_A,
\mathcal{I}_1),\mathcal{I}_2)\big)\E(\mathbf{1}_A,\mathcal{I}_2)d\mu \\
& \text{(and as $\E(\mathbf{1}_A, \mathcal{I}_1)(x)\leq 1$ we have
 $\E(\mathbf{1}_A, \mathcal{I}_2)\geq \E\big(\mathbf{1}_A\E(\mathbf{1}_A,
\mathcal{I}_1), \mathcal{I}_2\big)$} \\
&\geq \int \bigg(\E\big(\mathbf{1}_A\E(\mathbf{1}_A,
\mathcal{I}_1),\mathcal{I}_2)\big)\bigg)^2d\mu \\
&\geq \bigg(\int \mathbf{1}_A\E(\mathbf{1}_A,
\mathcal{I}_1d\mu\bigg)^2\geq \bigg(\int \big(\E(\mathbf{1}_A,
\mathcal{I}_1)\big)^2d\mu\bigg)^2\\
&\geq \bigg(\mathbf{1}_Ad\mu\bigg)^4= \mu(A)^4
\end{aligned}
\]
\end{proof}

The study of the case of three measure preserving transformations
seems much more complex.

\vskip1ex

\begin{lemma}
Let $(X, \mathcal{B}, \mu, T)$ be an invertible measure preserving
system on a finite measure space, $\mathcal{K}$ the $\sigma$
algebra spanned by the eigenfunctions of $T$ and $f$ a bounded
function. Let us denote by $X_f$ the set of full measure given by
the Wiener Wintner ergodic theorem such that for each $x\in X_f$
the averages
$$\frac{1}{N}\sum_{n=1}^N f(T^nx) e^{-2\pi i nt}$$ converge for
each t. For each $t\in \R$ let us denote by $E_t(f)$ the limit
function of these averages.
\begin{enumerate}
\item $E_t(f)$ is the projection of the function $f$ onto the
eigenspace of $T$ corresponding to the eigenvalue $e^{2\pi it}.$
In particular $E_0(f)$ is equal to $\E(f, \mathcal{I})$ the
conditional expectation with respect to the $\sigma$ algebra of
invariant sets for $T.$ \item If $t\neq s$ we have $\int
E_t(f)\overline{E_s(f)}d\mu = 0.$ \item If $e^{2\pi i\theta_k}$ is
the countable sequence of eigenvalues for $T$ and $e^{2\pi it_k}$
any countable set of distinct complex numbers then
$$\sum_{k=0}^{\infty}\|E_{t_k}(f)\|_2^2\leq
\sum_{k=0}^{\infty}\|E_{\theta_k}(f)\|_2^2\leq \|E(f,
\mathcal{K})\|_2^2$$
\end{enumerate}
\end{lemma}
\begin{proof}
 This is a simple consequence of the spectral theorem. If we denote by $P_t(f)$ the projection onto
 the eigenspace corresponding to the eigenvalue $e^{2\pi it}$ and by $\sigma_{f-P_t(f)}$ the spectral measure of the function
 $f-P_t(f)$ then we have
 \[\begin{aligned}
 &\frac{1}{N}\sum_{n=1}^N f(T^nx) e^{-2\pi i nt}=
 \frac{1}{N}\sum_{n=1}^N P_t(f)(T^nx) e^{-2\pi i nt}
 + \frac{1}{N}\sum_{n=1}^N [f-P_t(f)](T^nx) e^{-2\pi i nt} \\
 &= P_t(f)(x) + \frac{1}{N}\sum_{n=1}^N [f-P_t(f)](T^nx) e^{-2\pi i nt}
 \end{aligned}
 \]
As $$\lim_N \|\frac{1}{N}\sum_{n=1}^N [f-P_t(f)](T^nx) e^{-2\pi i
nt}\|_2^2 = \int \big|\frac{1}{N}\sum_{n=1}^N e^{2\pi
in(\theta-t)}\big|^2d\sigma_{f-P_t(f)}(\theta) =
\sigma_{f-P_t(f)}(\{0\}) = 0$$ we can conclude that $P_t(f)=
E_t(f).$ \vskip1ex

From this identification the remaining parts of the lemma follow
without difficulty. For the last part of the lemma we just need to
observe that $E_t(f) = 0 $ if $e^{2\pi it}$ is not an eigenvalue
of $T.$
\end{proof}

\vskip1ex

\noindent{\bf Remark}
 It is worth noticing that there is a key difference at the
 pointwise  level between $E_t(f)(x)$ and $P_t(f)(x).$ This
 difference highlights the difficulty one faces when dealing with
 ergodic versus not necessarily ergodic transformations.
The function $E_t(f)$ is defined off a single set of measure zero
for ALL $t\in \R.$ For each $t\in \R$ it is almost everywhere
equal to the function $P_t(f)(x)$ and so for each $t$ the $L^2$
functions $P_t(f)$ and $E_t(f)$ are equal. However we can not
claim that there is a universal null set off which one could write
that $E_t(f)(x) = P_t(f)(x)$ for all $t\in \R.$ One can look at
the example given in Proposition 7 below.
\begin{proposition}
Let $(X, \mathcal{B}, \mu, T_i)$, $1\leq i\leq 3,$ be three
measure preserving systems on the same finite measure space.There
exists a constant $0<\delta<1$ (independent from the $T_i$) such
that for all measurable set $A$ with measure $\mu(A)> 1-\delta$ we
have
$$\lim_N\frac{1}{N^2}\sum_{n,m=1}^N \mu(A\cap T_1^{-n}A\cap
T_2^{-m}A\cap T_3^{-n-m}A)\geq \frac{1}{2}\mu (A)^8$$
\end{proposition}
\begin{proof}
  We consider three measure preserving transformations $T_j,$ $1\leq
  j\leq 3$ and a measurable set $A$ with positive measure. We list the following
  notations and properties.
  \begin{enumerate}
  \item We denote by $E_{t}^j(\mathbf{1}_A)(x)$ the limit of
  $$\frac{1}{N}\sum_{n=1}^N \mathbf{1}_A(T_j^nx)e^{-2\pi int}$$
  for all $t\in \R$ off a single set of measure zero.
  \item For each $1\leq j \leq 3$ we consider the universal sets
  $X_{\mathbf{1}_A}^j$ such that $E_{t}^j(\mathbf{1}_A)(x)$ exists
  for all $t\in \R.$
  \item We consider an ergodic decomposition of $(X, \mathcal{B},
  \mu, T_3)$ with the measures $\mu_c$ where $d\mu = d\mu_cdc.$
  \item  We call $\mathcal{K}_c$ the Kronecker factor of $T_3$
  relative to the measure space $(X, \mathcal{B}, \mu_c).$
  The basis of eigenfunctions of $T_3$ relative to $\mu_c$ is denoted by
  $e_{k, c}.$ The constant function $1$ corresponds to $e_{0,c}.$
  \item The eigenvalue corresponding to the eigenfunction $e_{k, c}$
  is $e^{-2\pi i\theta_{k,c}}.$
  \item By Birkhoff pointwise ergodic theorem combined with the
  disintegration of $\mu$ we have for a.e. $c$ for $\mu_c$ a.e. $y$
   $\dis \lim_N \frac{1}{N}\sum_{n=1}^N \mathbf{1}_A(T_3^ny) =
   E(\mathbf{1}_A,\mathcal{I})(y)= E_0^3(\mathbf{1}_A)(y),$
   where $\mathcal{I}$ denotes the $\sigma$ algebra of invariant
   sets with respect to $\mu.$
  \end{enumerate}
   As a consequence of the ergodicity of $T_3$ with respect to $\mu_c$ we have
   $$\E(\mathbf{1}_A, \mathcal{K}_c)= \sum_{k=0}^{\infty}\big(\int
   \mathbf{1}_A\overline{e_{k, c}}d\mu_c \big)e_{k,c}.$$
  We disintegrate the measurable sets $X_{\mathbf{1}_A}^j$ with respect to the
  measure $d\mu_c.$ We obtain for $\mu_c$ a.e. $y$ for all $t\in \R$
  the pointwise convergence of the averages
  $$\frac{1}{N}\sum_{n=1}^N \mathbf{1}_A(T_j^ny)e^{-2\pi int}.$$
  This is crucial for our method as with respect to the measure $\mu_c$
  the transformations $T_1$ and $T_2$ are not necessarily measure
  preserving.
  \vskip1ex
  For each eigenfunction $e_{k,c}$ we have

  \[\begin{aligned}
  &\lim_N\frac{1}{N^2}\sum_{n,m=1}^N
  \mathbf{1}_A(T_1^ny)\mathbf{1}_A(T_2^my)e_{k,c}(T_3^{n+m}y)\\
  &= e_{k,c}(y)\lim_N\frac{1}{N^2}\sum_{n,m=1}^N
  \mathbf{1}_A(T_1^ny)e^{-2\pi
  in\theta_{k,c}}\mathbf{1}_A(T_2^my)e^{-2\pi im\theta_k}
  = e_{k,c}(y)E^1_{\theta_{k,c}}(\mathbf{1}_A)(y)E^2_{\theta_{k,c}}(\mathbf{1}_A)(y)
   \end{aligned}
   \]

  As a consequence we have
  \[\begin{aligned}
  &\lim_N \frac{1}{N^2}\sum_{n,m=1}^N \mu_c\big(A\cap T_1^nA\cap T_2^mA\cap T_3^{n+m}A\big)
  =\sum_{k=0}^{\infty} \big(\int \mathbf{1}_A
\overline{e_{k,c}}d\mu_c\big)
  \big(\int \mathbf{1}_A
  e_{k,c}E^1_{\theta_{k,c}}(\mathbf{1}_A)E^2_{\theta_{k,c}}(\mathbf{1}_A)d\mu_c\big)\\
  &= \mu_c(A)\big(\int \mathbf{1}_A
  E^1_{0}(\mathbf{1}_A)E^2_{0}(\mathbf{1}_A)d\mu_c\big)+ \sum_{k=1}^{\infty} \big(\int \mathbf{1}_A
\overline{e_{k,c}}d\mu_c\big)
  \big(\int \mathbf{1}_A
  e_{k,c}E^1_{\theta_{k,c}}(\mathbf{1}_A)E^2_{\theta_{k,c}}(\mathbf{1}_A)d\mu_c\big)\\
  \end{aligned}
  \]
 The first term of the previous line is for a.e. $c$ equal to
 $$ \E( \mathbf{1}_A, \mathcal{I}_3)\int \mathbf{1}_A(y)\E( \mathbf{1}_A,
 \mathcal{I}_1)\E( \mathbf{1}_A, \mathcal{I}_2)d\mu_c.$$
 In view of the constance of $\E( \mathbf{1}_A, \mathcal{I}_3)$ with
 respect to $\mu_c,$
 this last term can be written as
 $$\int \E( \mathbf{1}_A, \mathcal{I}_3)\mathbf{1}_A(y)\E( \mathbf{1}_A,
 \mathcal{I}_1)\E( \mathbf{1}_A, \mathcal{I}_2)d\mu_c.$$
  Integrating with respect to $dc$ and using properties of the conditional expectation we get
 \[\begin{aligned}
 &\int \E( \mathbf{1}_A, \mathcal{I}_3)\mathbf{1}_A(y)\E( \mathbf{1}_A,
 \mathcal{I}_1)\E( \mathbf{1}_A, \mathcal{I}_2)d\mu_cdc = \int \E( \mathbf{1}_A, \mathcal{I}_3)\mathbf{1}_A\E( \mathbf{1}_A,
 \mathcal{I}_1)\E( \mathbf{1}_A, \mathcal{I}_2)d\mu(x)\\
 &\geq \int \E\big(\mathbf{1}_A\E( \mathbf{1}_A, \mathcal{I}_3)\E( \mathbf{1}_A,
 \mathcal{I}_2), \mathcal{I}_1\big)\mathbf{1}_A d\mu = \int \E\big(\mathbf{1}_A\E( \mathbf{1}_A, \mathcal{I}_3)\E( \mathbf{1}_A,
 \mathcal{I}_2), \mathcal{I}_1\big)\E(\mathbf{1}_A, \mathcal{I}_1)d\mu\\
 &\geq \big(\int \big(\E\big(\mathbf{1}_A\E( \mathbf{1}_A, \mathcal{I}_3)\E( \mathbf{1}_A,
 \mathcal{I}_2), \mathcal{I}_1\big)\big)^2 d\mu\geq \bigg(\int \mathbf{1}_A\E( \mathbf{1}_A, \mathcal{I}_3)\E( \mathbf{1}_A,
 \mathcal{I}_2)d\mu\bigg)^2 \\
 &\geq \bigg(\int \E\big(\mathbf{1}_A,\E( \mathbf{1}_A,
 \mathcal{I}_3), \mathcal{I}_2\big)\E(\mathbf{1}_A, \mathcal{I}_2)
 d\mu\bigg)^2 \\
 &= \bigg( \int \big(\E\big(\mathbf{1}_A,\E( \mathbf{1}_A,
 \mathcal{I}_3), \mathcal{I}_2\big)\big)^2d\mu\bigg)^2\geq \bigg(\int \mathbf{1}_A \E( \mathbf{1}_A,
 \mathcal{I}_3)d\mu\bigg)^4 \geq \mu(A)^8.
 \end{aligned}
 \]
  This shows that after integration with respect to $dc$  $\dis \mu_c(A)\big(\int \mathbf{1}_A
  E^1_{0}(\mathbf{1}_A)E^2_{0}(\mathbf{1}_A)d\mu_c\big)$ is bounded below by $\mu(A)^8.$
 Our proof will be complete if one can show that if $\mu(A) > 1- \delta$ for some universal $0<\delta<1$
 then
 \begin{equation}
 \int |(I)_c|dc = \int \bigg|\int \sum_{k=1}^{\infty} \big(\int \mathbf{1}_A
\overline{e_{k,c}}d\mu_c\big)
  \big(\int \mathbf{1}_A
  e_{k,c}E^1_{\theta_{k,c}}(\mathbf{1}_A)E^2_{\theta_{k,c}}(\mathbf{1}_A)d\mu_c\big)\bigg|dc \leq
  \frac{1}{2}\mu(A)^8.
  \end{equation}
By Cauchy-Schwartz's inequality we have
$$|(I)_c|\leq \bigg(\sum_{k=1}^{\infty} \big|\int
\mathbf{1}_A\overline{e_{k,c}}d\mu_c\big|^2\bigg)^{1/2}\bigg(\sum_{k=1}^{\infty}\big|\int\mathbf{1}_A
  e_{k,c}E^1_{\theta_{k,c}}(\mathbf{1}_A)E^2_{\theta_{k,c}}(\mathbf{1}_A)d\mu_c\big|^2\bigg)^{1/2}.$$
   The vectors $e_{k,c}$ form an orthonormal basis of $L^2(X,
  \mathcal{K}_c, \mu_c)$ because $T_3$ on this space is ergodic. Thus we have
\begin{equation}
\bigg(\sum_{k=1}^{\infty} \big|\int
\mathbf{1}_A\overline{e_{k,c}}d\mu_c\big|^2\bigg)^{1/2}=
\bigg(\int \big|\E(\mathbf{1}_A, \mathcal{K}_c)\big|^2d\mu_c -
\mu_c(A)^2\bigg)^{1/2}\leq \big(\mu_c(A)- \mu_c(A)^2\big)^{1/2}
\end{equation}
The second term $\dis (II)_c =
\bigg(\sum_{k=1}^{\infty}\big|\int\mathbf{1}_A
  e_{k,c}E^1_{\theta_{k,c}}(\mathbf{1}_A)E^2_{\theta_{k,c}}(\mathbf{1}_A)d\mu_c\big|^2\bigg)^{1/2}$
  can also be bounded above by
  $$\bigg(\sum_{k=1}^{\infty} \bigg(\int
  \mathbf{1}_A\big|E^1_{\theta_{k,c}}(\mathbf{1}_A)\big|^2d\mu_c\bigg)^2\bigg)^{1/2}\bigg(\sum_{k=1}^{\infty} \bigg(\int
  \mathbf{1}_A\big|E^2_{\theta_{k,c}}(\mathbf{1}_A)\big|^2d\mu_c\bigg)^2\bigg)^{1/2}.$$
 Using Lemma 9 part (3), this last term is bounded above by
 $$\bigg( \int |\E(\mathbf{1}_A, \mathcal{K}_c)|^2d\mu_c -\mu_c(A)^2\bigg)^{1/2}
 \bigg( \int |\E(\mathbf{1}_A, \mathcal{K}_c)|^2d\mu_c
 -\mu_c(A)^2\bigg)^{1/2}$$ which is equal to
\begin{equation}
 \bigg( \int |\E(\mathbf{1}_A, \mathcal{K}_c)|^2d\mu_c
 -\mu_c(A)^2\bigg)\leq (\mu_c(A) - \mu_c(A)^2)
 \end{equation}
 Combining the bounds found in (2) and in (3) we get

 $\dis |(I)_c|\leq (\mu_c(A) - \mu_c(A)^2)^{3/2}.$
 As $\dis (\mu_c(A) - \mu_c(A)^2)^{3/2}\leq (\mu_c(A) - \mu_c(A)^2)$ and $\int \mu_c(A)^2dc \geq (\int \mu_c(A)dc )^2,$
 integrating with respect to c we obtain
 $$\int |(I)_c| dc\leq \int (\mu_c(A) - \mu_c(A)^2)^{3/2}dc\leq \int (\mu_c(A) - \mu_c(A)^2)dc \leq \mu(A) - \mu(A)^2.$$
 Going back to (1) we will reach our conclusion if we can find
 $0<\delta<1$ such that
 $$\mu(A) - \mu(A)^2\leq \frac{1}{2} \mu(A)^8,$$
 for all measurable set $A$ with measure greater or equal to
 $1- \delta.$
 This is an easy consequence of the uniqueness of the root for the
 polynomial $\dis 1/2 x^7 + x -1$ on $(0,1).$

\noindent{\bf Remark} The constant $\frac{1}{2}$ for the lower
bound $\frac{1}{2}\mu(A)^{8}$ is certainly not optimal. Following
the same path one can show that
$$\lim_N\frac{1}{N^2}\sum_{n,m=1}^N \mu(A\cap T_1^{-n}A\cap
T_2^{-m}A\cap T_3^{-n-m}A)>0$$ for all measurable set $A$ when
$\mu(A)>\beta$ where $\beta$ is the root of $x^7 +x -1$ on
$(0,1).$

  \end{proof}
  \section{On the almost everywhere convergence of weighted averages}
  \subsection{The averages $\frac{1}{N^2}\sum_{n,m=1}^Na_nb_mf(T^{n+m}x)$}

  \vskip1ex

  Our goal is to prove first that the ergodicity assumptions are
  necessary in Lemma 3. We recall that we
  denote by $\mathcal{K}$ the $\sigma$ algebra generated by the
  eigenfunctions of a measure preserving transformation. Even
  without the ergodicity assumption this $\sigma$-algebra is well
  defined.
 \begin{proposition}
 There exists a non ergodic measure preserving system $(Y,\mathcal{B}, \nu,
 S),$ a function $f\in L^{\infty}(\nu)\cap \mathcal{K}^{\perp}$ such that
 for $\nu$ a.e. $y$ we can find bounded sequences $a_n$ and $b_n$ such that
the averages $$\frac{1}{N^2}\sum_{n, m=1}^N a_nb_mf(S^{n+m}y)$$ do
not converge when $N$ tends to $\infty$. In other words Lemma 3 is
false if we remove the ergodicity assumption.
\end{proposition}
\begin{proof}
  Let $S(x, y) = (x+ \al, x+y)$ be the ergodic measure preserving
  transformation defined on the two Torus where $\al$ is an irrational number.
   We consider the measure preserving
  transformation $T = S\times S$ on $\T^4$ defined as
  $$T(x_1, x_2, x_3, x_4) = (x_1 + \al, x_1 + x_2, x_3 +\al, x_3 +
  x_4).$$
  The transformation $T$ is not ergodic and the Kronecker factor
  ($\sigma$ algebra spanned by the eigenfunctions of $T$)
  corresponds to the functions depending on the first and third
  coordinates $x_1$ and $x_3$. This is because the eigenfunctions
  of $S$ depend on their first coordinates (see also Lemma 4.18 in
  \cite{Furstenberg} on the way in general the eigenfunctions of $T$ are
  created from those of $S$)
  Consider the function $f(x_1, x_2, x_3, x_4)= e^{-2\pi i x_2}e^{2\pi i x_4}.$
  This function belongs to $\mathcal{K}^{\perp}.$
  We have $f(T^{n+m}(x_1, x_2, x_3, x_4))= e^{2\pi i(x_4- x_2 +
  (n+m)(x_3-x_1)}.$  Let us assume that Lemma 3 was true without ergodicity
  assumption then we could find a set of
  full measure such that for a.e. $x_1, x_2, x_3, x_4$ in this set
  and for all bounded sequences $a_n$ and $b_n$ we would have
  $$\lim_N\frac{1}{N^2}\sum_{n, m=0}^{N-1}a_nb_mf(T^{n+m}(x_1, x_2, x_3,
  x_4))=\lim_N\frac{1}{N^2}\sum_{n, m=0}^{N-1}a_nb_me^{2\pi i(x_4- x_2 +
  (n+m)(x_3-x_1)}= 0.$$
  To disprove this we can take a bounded sequence $v_n$ such that
  the averages $\frac{1}{N}\sum_{n=0}^{N-1}v_n$ diverge. Then we
  can take $a_n= v_ne^{-2\pi in(x_3- x_1)},$ and
  $b_m= e^{-2\pi im(x_3- x_1)}.$ As $$\frac{1}{N^2}\sum_{n, m=0}^{N-1}a_nb_mf(T^{n+m}(x_1, x_2, x_3,
  x_4))= \frac{1}{N}\sum_{n=0}^{N-1}v_n e^{2\pi i(x_4- x_2 )},$$
 this shows that Lemma 3 is false once we remove the ergodicity
  assumption. This ends the proof of Proposition 7.
  \end{proof}
\vskip1ex \noindent{\bf{Remarks 1}}
\begin{enumerate}
\item
  Proposition 7 shows that Lemma 3 as stated is quite sharp as one can not even
  expect to have the convergence of the averages
  $\dis \frac{1}{N^2}\sum_{n,m=0}^{N-1}a_nb_mf(T^{n+m}x)$ as in this
  example they are equal to $\dis \frac{1}{N}\sum_{n=1}^N v_n.e^{2\pi i(x_4- x_2)}$
\item The same measure preserving system can be used to show that
the uniform Wiener Wintner ergodic theorem is no longer valid if
$T$ is not ergodic. By this we mean that if we denote by
$\mathcal{K}$ the $\sigma$ algebra
   spanned by the eigenfunctions of $T$ then we do not have in general for
  functions $f\in \mathcal{K}^{\perp},$
  $$\lim_N\sup_t |\frac{1}{N}\sum_{n=0}^{N-1} f(T^nx) e^{2\pi
  int}|= 0.$$
\item
  As indicated earlier the norm convergence holds without difficulty as the next proposition shows.
  We give the proof just for the sake of completeness and to show
 the difference between the pointwise and norm convergence.
\end{enumerate}

  \vskip1ex

 \begin{definition}
 We will denote by $\mathcal{WW}_1$ the set of bounded sequences $a=(a_n)$ of
 scalars such that $\lim_N\frac{1}{N}\sum_{n=1}^Na_ne^{2\pi int}$
 exists for each $t\in \R.$
 \end{definition}
\begin{proposition}
Let $T$ be a unitary operator and let $a= (a_n)$ and $b=(b_m)$ be
bounded sequences. Then the averages
$$\frac{1}{N^2}\sum_{n, m=1}^N a_nb_mT^{n+m}$$ converge in norm
 if a and b belong to $\mathcal{WW}_1$
\end{proposition}
\begin{proof}
It is a simple consequence of the spectral theorem. If we denote
by $\sigma_f$ the spectral measure of the function $f$ with
respect to $T$ then we have
\[\begin{aligned}
&\|\frac{1}{N^2}\sum_{n, m=1}^N a_nb_mT^{n+m}-
\frac{1}{M^2}\sum_{n, m=1}^M a_nb_mT^{n+m}\|^2\\
&=\int \big|\frac{1}{N^2}\sum_{n, m=1}^N a_nb_me^{2\pi i(n+m)t}-
\frac{1}{M^2}\sum_{n, m=1}^M a_nb_me^{2\pi
i(n+m)t}\big|^2d\sigma_f(t)\\
&= \int \big|\frac{1}{N}\sum_{n=1}^N a_n e^{2\pi
int}\frac{1}{N}\sum_{m=1}^Nb_m e^{2\pi i mt}-
\frac{1}{N}\sum_{n=1}^N a_n e^{2\pi int}\frac{1}{N}\sum_{m=1}^Nb_m
e^{2\pi i mt}\big|^2d\sigma_f(t)
\end{aligned}
\]
which easily shows that the averages form a Cauchy sequence.
\end{proof}

\subsection{Higher order averages}
\vskip1ex

Proposition 2 shows that if the transformation $T$ is ergodic and
the function $f\in L^2$ then for $\mu$ a.e. $x$ the averages
$$\frac{1}{N^2}\sum_{n,m=0}^{N-1} a_nb_mf(T^{n+m}x)$$ converge for all sequences $a= (a_n), b = (b_n)$ that belong to
 $\mathcal{WW}_1$.
\vskip1ex The next proposition shows that the class
$\mathcal{WW}_1$ does not characterize those bounded sequences for
which the similar averages for seven terms converge a.e. even
under the condition of ergodicity of the transformation
 \vskip1ex
\begin{proposition}
There exists an ergodic dynamical system $(X, \mathcal{A}, \mu,
T)$ and a function $f\in L^{\infty}(\mu)$ such that for $\mu$ a.e.
$x$ we can find bounded sequences $A_i= (a_{n,i})\in
\mathcal{WW}_1$ for whcih the averages
$$M_N(A_1, A_2,...A_6,f)(x)= \frac{1}{N^3}\sum_{n,
m,p=0}^{N-1}a_{1,p}a_{2,n}a_{3,p+n}a_{4,m}a_{5,n+m}a_{6,
p+m}f(T^{n+m+p}x)$$ do not converge.
\end{proposition}
\begin{proof}

We consider the sequence $v_n$ with values 1 or -1 such that the
averages $\frac{1}{N}\sum_{n=1}^N v_n$ diverge. The sequence is
built from longer and longer stretches of 1 and -1 so that the
averages get close to 1 then close to -1 and so on. We extend
$v_n$ to negative indices by putting $v_{-n}= v_n$. We can observe
that this sequence has a correlation in the sense that for any
$h\in \Z$ the averages $\dis \frac{1}{N}\sum_{n=1}^N v_n
\overline{v_{n+h}}$ converge to a scalar $\gamma (h)$. Simple
considerations show that the limit for all $h$ is equal to one.
The quantity $\gamma(h)$ represents the $h$ Fourier coefficients
of a positive measure $\sigma$ that is equal then to the Dirac
measure at zero, $\delta_0$, a discrete measure.

\vskip1ex

  We take now an irrational number $\al$ . We claim that the
sequence $a_n= v_ne^{2\pi in^2\al}$ belongs to $\mathcal{WW}_1.$
To see this first one can observe that the sequence $e^{2\pi
in^2\al}e^{2\pi int}$ does have a correlation; for each $h\in \Z$
the limit of $\dis\frac{1}{N}\sum_{n=1}^N e^{2\pi i[(n^2-
(n+h)^2]\al}e^{2\pi nt -(n+h)t}$ is equal to zero for $h\ne 0$ and
zero otherwise. Therefore the measure associated with these
Fourier coefficients is Lebesgue measure, $m$. As a consequence of
the Affinity principle the measures $m$ and $\delta_0$ being
orthogonal we have for each $t\in \R$
$$\lim_N \frac{1}{N}\sum_{n=1}^N v_ne^{2\pi in^2\al}e^{2\pi int}=
0.$$ Thus we have shown that the sequence $a_n= v_ne^{2\pi
in^2\al}$ belongs to $\mathcal{WW}_1.$

\vskip1ex

We consider the ergodic measure preserving transformation $S(x, y)
= (x+ \al,x+y)$ defined on the two Torus where $\al$ is the
irrational number used to define the sequence $a_n.$ Our goal is
to prove that for the function $f(x,y) = e^{4\pi iy}$ it is
impossible to find a set of full measure off which for all six
bounded sequences $A_i= (a_{i,n}),$ $1\leq i\leq 6$ the averages
$$\frac{1}{N^3}\sum_{n,
m,p=0}^{N-1}a_{1,p}a_{2,n}a_{3,p+n}a_{4,m}a_{5,n+m}a_{6,
p+m}f(S^{n+m+p}(x,y))$$ converge. To reach this conclusion we can
use the simple equality
$$(n+m+p)^2 = (n+m)^2 + (n+p)^2 + (m+p)^2 - n^2 -p^2 -m^2.$$
We have $f(T^{n+m+p}(x,y))= e^{4\pi i(y +
(n+m+p)(x-\al/2)}.e^{2\pi
 i(n+m+p)^2\al}.$ As a consequence if we take $a_{1,p}= e^{2\pi
 ip^2\al},$ $a_{2,n}= v_n e^{2\pi in^2\al},$ $a_{3, p+n}= e^{-2\pi
 i(p+n)^2\al}e^{-2\pi i(p+n)(x-\al/2)},$ $a_{4,m}= e^{2\pi
 im^2\al},$ $a_{5,n+m}= e^{-2\pi i(n+m)^2\al}e^{-2\pi
 i(n+m)(x-\al/2)},$ and $a_{6,p+m}= e^{-2\pi i(p+m)^2\al}e^{-2\pi
 i(p+m)(x-\al/2)},$ then
$$\frac{1}{N^3}\sum_{n,
m,p=0}^{N-1}a_{1,p}a_{2,n}a_{3,p+n}a_{4,m}a_{5,n+m}a_{6,
p+m}f(S^{n+m+p}(x,y))= \frac{1}{N}\sum_{n=0}^{N-1} v_ne^{4\pi
iy}.$$ Therefore the averages
$$\frac{1}{N^3}\sum_{n,
m,p=0}^{N-1}a_{1,p}a_{2,n}a_{3,p+n}a_{4,m}a_{5,n+m}a_{6,
p+m}f(S^{n+m+p}(x,y))$$ do not converge. (The arguments in the
previous paragraphs can also be used to show that each sequence
$A_i\in \mathcal{WW}_1.$)
\end{proof}

 \end{document}